\documentclass[american, 5p]{elsarticle}
\pdfoutput=1
\journal{}
\biboptions{longnamesfirst, sort&compress}
\makeatletter
\def\ps@pprintTitle{%
   \def\@oddfoot{\reset@font\phantom{\today}\hfil\thepage\hfil\today}
}
\makeatother

\usepackage[utf8]{inputenc}
\usepackage[T1]{fontenc}
\usepackage{lmodern}
\usepackage[babel=true]{microtype}

\usepackage{babel}
\usepackage{amsfonts, amsmath, amsopn, amssymb, bm, braket}
\usepackage{enumitem}
\usepackage[caption=false]{subfig}
\usepackage[export]{adjustbox}

\usepackage{pgfplots, tkz-euclide}
\pgfplotsset{compat=newest}
\usetikzlibrary{calc}
\usepackage{siunitx}
\usepackage{hyperref}
\hypersetup{%
   pdfauthor  = {Lennard Kamenski},
   pdftitle   = {Sharp bounds on the smallest eigenvalue of finite element equations
   with arbitrary meshes without regularity assumptions},
   pdfsubject = {Preprint}
}

\DeclareMathOperator{\supp}{supp}
\newcommand{\abs}[1]{\lvert#1\rvert}
\newcommand{\Abs}[1]{\left\lvert#1\right\rvert}
\newcommand{\bigAbs}[1]{\big\lvert#1\big\rvert}

\newcommand{\Norm}[1]{\left\lVert#1\right\rVert}
\newcommand{\bx}{\bm{x}}
\newcommand{\bxi}{\bm{\xi}}
\newcommand{\bu}{\bm{u}}
\newcommand{\dx}{\,\mathrm{d}\bm{x}}
\newcommand{\dxi}{\,\mathrm{d}\bm{\xi}}
\newcommand{\D}{\mathbb{D}}
\newcommand{\cN}{\mathcal{N}}

\newcommand\Th{\mathcal{T}_h}

\usepackage[capitalize]{cleveref} % For smart references
\crefname{equation}{}{}
\crefname{figure}{Fig.}{Figs.}
\Crefname{figure}{Figure}{Figures}
\crefname{section}{section}{sections}
\Crefname{section}{Section}{Sections}
\crefname{assumption}{Assumption}{Assumptions}

\usepackage{amsthm}
\newtheorem{theorem}{Theorem}
\newtheorem{assumption}{Assumption}
\theoremstyle{definition}
\newtheorem{remark}{Remark}
\newtheorem{meshexample}{Mesh example}

\newcommand{\PlotMesh}[2]{%
  \begin{tikzpicture}
    \begin{axis}[%
      width = 1.0\linewidth, scale=#2, scale only axis,
      trim axis left, trim axis right, axis equal image,
      axis x line*={left}, axis y line*={left},
      xmin={-0.03}, xmax={1.03},
      xtick={{0.0,1.0}}, xticklabels={{$0$,$1$}},
      xticklabel style={font={{\fontsize{8 pt}{11 pt}\selectfont}}},
      ymin={-0.03}, ymax={1.03},
      ytick={{0.0,1.0}}, yticklabels={{$0$,$1$}},
      yticklabel style={font={{\fontsize{8 pt}{11 pt}\selectfont}}},
      ]
      \addplot[patch, mesh, patch type=triangle, color=black,
        patch table = {#1.ele},
      ]
      table [x index=0, y index=1] {#1.coor};
      \end{axis}
  \end{tikzpicture}%
}

\newcommand{\PlotResults}[4][]{%
  \begin{tikzpicture}
    \begin{loglogaxis}[%
      width={1.0\linewidth}, height={0.7\linewidth}, scale=#4,
      grid=major, trim axis left, trim axis right,
      legend columns=2,
      legend cell align={left},
      legend style={%
        font={{\fontsize{8 pt}{10 pt}\selectfont}},
        at={(0.02, 0.02)}, anchor={south west}%
      },
      anchor={north west}, xshift={1.0mm}, yshift={-1.0mm},
      xlabel={#3},
      xlabel style={
        font={{\fontsize{8 pt}{10 pt}\selectfont}},
       },
      xticklabel style={
        font={{\fontsize{8 pt}{10 pt}\selectfont}},
      },
      x grid style={%
        line width={0.5},
      },
      axis x line*={left},
      ylabel style={%
        font={{\fontsize{8 pt}{10 pt}\selectfont}},
      },
      yticklabel style={
        font={{\fontsize{8 pt}{10 pt}\selectfont}},
      },
      y grid style={ line width={0.5},
      },
      axis y line*={left},
      ]
      \addplot[%
        color={rgb,1:red,0.0;green,0.6056;blue,0.9787},
        line width={1}, solid, mark={*}, mark size={2.25 pt},
        mark options={%
          color={rgb,1:red,0.0;green,0.6056;blue,0.9787},
          fill={rgb,1:red,0.0;green,0.6056;blue,0.9787},
        }
      ]
      table [x index=0, y=lmin] {#2};
      \addlegendentry {$\lambda_{\min}$}
      \addplot[%
        color={rgb,1:red,0.8889;green,0.4356;blue,0.2781},
        line width={1}, dashed, mark={star}, mark size={3.0 pt},
        mark options={%
          color={rgb,1:red,0.8889;green,0.4356;blue,0.2781},
          fill={rgb,1:red,0.8889;green,0.4356;blue,0.2781},
        }
      ]
      table [x index=0, y=lminE] {#2};
      \addlegendentry {$\bar{\lambda}$}
      \addplot[%
        color={rgb,1:red,0.2422;green,0.6433;blue,0.3044},
        line width={1}, dashdotted, mark={triangle*}, mark size={3.0 pt},
        mark options={%
          color={rgb,1:red,0.2422;green,0.6433;blue,0.3044},
          fill={rgb,1:red,0.2422;green,0.6433;blue,0.3044},
          solid,
        }
      ]
      table [x index=0, y=lminKHX] {#2};
      \addlegendentry {$\bar{\lambda}_{KHX}$}
      \addplot[%
        color={rgb,1:red,0.7644;green,0.4441;blue,0.8243},
        line width={1}, dotted, mark={diamond*}, mark size={3.0 pt},
        mark options={%
          color={rgb,1:red,0.7644;green,0.4441;blue,0.8243},
          fill={rgb,1:red,0.7644;green,0.4441;blue,0.8243},
          solid,
        }
      ]
      table [x index=0, y=lminGM] {#2};
      \addlegendentry {$\bar{\lambda}_{GM}$}
      #1
  \end{loglogaxis}
\end{tikzpicture}

}

\begin{document}

\begin{frontmatter}
\title{%
   Sharp bounds on~the smallest eigenvalue of~finite element equations\\
   with~arbitrary meshes without regularity assumptions%
}

\author{Lennard Kamenski}

%--- abstract --------------------------------------------------------
\begin{abstract}
A proof for the lower bound is provided for the smallest eigenvalue of finite element equations with arbitrary conforming simplicial meshes.
The bound has a similar form to the one by Graham and McLean [\emph{SIAM J.\ Numer.\ Anal.}, 44 (2006), pp.~1487--1513] but doesn't require any mesh regularity assumptions, neither global nor local.
In particular, it is valid for highly adaptive, anisotropic, or nonregular meshes without any restrictions.
In three and more dimensions, the bound depends only on the number of degrees of freedom $N$ and the Hölder mean $M_{1-d/2} (\lvert \tilde{\omega} \rvert / \lvert \omega_i \lvert)$ taken to the power $1-2/d$, $\lvert \tilde{\omega} \rvert$ and $\lvert \omega_i \rvert$ denoting the average mesh patch volume and the volume of the patch corresponding to the $i^{\text{th}}$ mesh node, respectively.
In two dimensions, the bound depends on the number of degrees of freedom $N$ and the logarithmic term $(1 + \lvert \ln (N \lvert \omega_{\min} \rvert) \rvert)$, $\lvert \omega_{\min} \rvert$ denoting the volume of the smallest patch.
Provided numerical examples demonstrate that the bound is more accurate and less dependent on the mesh nonuniformity than the previously available bounds.
\end{abstract}

%--- keywords --------------------------------------------------------
\begin{keyword}
  finite element method, simplicial meshes, conditioning, eigenvalue estimates, extreme eigenvalues
  \MSC[2020] 65N30, 65N50, 65N22, 65F15%

  \vspace{0.6\baselineskip}%

  \begin{tabular}{@{}p{0.13\textwidth}@{}p{0.87\textwidth}@{}}
    \footnotesize{}%
    \includegraphics[width=0.115\textwidth, valign=t]{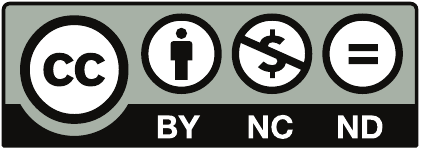}
  & \footnotesize{}%
    This is a preprint of a contibution published  in
    \emph{SIAM~J.~Numer.~Anal.}, 59(2) (2021), pp.~983--997.
   \vspace{0.3\baselineskip}
   \newline{}
   \copyright~2021. Licensed under CC-BY-NC-ND~4.0.
   The final version is available online at \url{https://doi.org/10.1137/19M128034X}.
   \end{tabular}%
\end{keyword}
\end{frontmatter}

%% section 1 %%%%%%%%%%%%%%%%%%%%%%%%%%%%%%%%%%%%%%%%%%%%%%%%%%%%%%%%%
\section{Introduction}

Highly adapted or misshaped meshes can have a strong influence on the conditioning of the finite element equations, which has an impact on convergence properties and accuracy of iterative methods for solving the resulting linear systems.
Thus, bounds on the condition number of finite element stiffness matrix and its extreme eigenvalues have been of theoretical and practical interest for a long time.
In most of the available work, some local or global mesh regularity properties are required in order to bound the constants during the derivation of the estimates, which causes degeneration of these estimates in the case of extreme mesh geometries or strong irregular anisotropic adaptation.
The aim of this work is to demonstrate that such mesh regularity assumptions are not necessarily required.

The estimation of the largest eigenvalue is well understood and it is easy to show that the largest eigenvalue is bounded by a multiple (with a constant depending on mesh connectivity) of the maximum of the largest eigenvalues of the local stiffness matrices, e.g.,~\cite{Fri73}.
Sharp bounds in terms of mesh geometry are available for both isotropic~\cite{AinMcLTra99, DuWanZhu09, GraMcL06, She02a} and anisotropic~\cite{HuaKamLan15,KamHuaXu14,She02a} diffusion.
Moreover, the largest diagonal element of the stiffness matrix is a good estimate for the largest eigenvalue: it is tight within a constant factor bounded by the maximal number of basis functions per element~\cite{HuaKamLan15}.
For the smallest eigenvalue, a good bound proved to be much more difficult to achieve.
For not quasi-uniform meshes, the bound on the smallest eigenvalue is often considered in terms of the volume of the smallest mesh element, e.g.,~\cite[Remark~9.11]{ErnGue04},~\cite{DuWanZhu09,Fri73,ZhuDu11, ZhuDu14}, which is clearly much too pessimistic.

A better bound was proven by Bank and Scott~\cite{BanSco89} for general isotropic meshes in $d \geq 2$ dimensions.
They showed that under a proper scaling the conditioning of finite element equations with local refinement does not degrade significantly and is comparable to that on a uniform mesh with a comparable number of degrees of freedom.
A more general result for elliptic bilinear forms on Sobolev spaces of real index $m \in [-1,1]$ with shape-regular meshes in $d \geq 2$ dimensions was derived by Ainsworth, McLean, and Tran~\cite{AinMcLTra99,AinMcLTra00}.

Their result was generalized by Graham and McLean~\cite{GraMcL06} to general meshes satisfying the following weak local regularity condition requiring the neighboring mesh elements to be comparable in size and shape.
\begin{assumption}[see \protect{\cite[Assumption~3.2]{GraMcL06}}]%
\label{ass:grmcl}
There exist positive constants $F$, $G$, $H$, and $M$ such that
\begin{align*}
  &\frac{h_K}{h_{K'}} \le F,
  ~ \frac{\rho_K}{\rho_{K'}} \le G,
  ~ \frac{\abs{K}}{\abs{K'}} \le H
  \quad \forall K, K' \colon \overline{K}\cap\overline{K'} \neq \emptyset,
  \\
  &\max_{i} \# \left\{K : \bm{x}_i \in \overline{K} \right\} \le M
  .
  %\label{eq:ass2}
\end{align*}
\end{assumption}
For the case of finite elements in $d \ge 3$ dimensions, the obtained bound for the smallest eigenvalue is
\begin{equation}
      \lambda_{\min}(A) \gtrsim
      {\left(  \sum\limits_{i\in\cN}  \Abs{\omega_i}^{1-\frac{d}{2}} \right)}^{-\frac{2}{d}}
      ,
      %\tag{GM}
   \label{eq:gramcl}
\end{equation}
$\Abs{\omega_i}$ denoting the volume of the patch corresponding to the $i^{\text{th}}$ mesh node~\cite[Lemma~5.2]{GraMcL06}.
However, looking at the constant in the bound reveals that the more correct form of this bound is (cf.~\cite[Eq.~(5.16) and the following]{GraMcL06})
\begin{equation}
      \lambda_{\min}(A) \gtrsim
      \frac{1}{{\left(M H\right)}^{\frac{d-2}{d}}}
      {\left(  \sum\limits_{i\in\cN}  \Abs{\omega_i}^{1-\frac{d}{2}} \right)}^{-\frac{2}{d}}
      .
   %\tag{GM}
   \label{eq:gramcl2}
\end{equation}
For meshes fulfilling \cref{ass:grmcl}, $M$ and $H$ remain bounded but for highly anisotropic meshes, e.g., Shishkin or Bakhvalov type, especially in higher dimensions, $H$ can become large, as shown in numerical experiments in \cref{sec:examples}.
Hence, for highly adaptive anisotropic meshes, $H$ may be considered as a mesh-dependent part of the bound.
$M$ usually remains bounded even for highly adaptive meshes, although it can grow in the case of a local mesh degeneracy commonly known as a \emph{spider}.%
\footnote{This is one central node connected to many neighbouring nodes, resembling a mesh of a spider.}

%--- notation figure---
\begin{figure}[t]%
   \centering{}%
   \begin{tikzpicture}[scale = 0.87]%
      \tikzstyle{every node}=[]%
      %
      % patch nodes
      %
      \path ( 0.0,  0.0) coordinate (N0);
      \path ( 2.0,  0.0) coordinate (N1);
      \path ( 0.6,  1.8) coordinate (N2);
      \path (-1.5,  1.3) coordinate (N3);
      \path (-1.4, -1.0) coordinate (N4);
      \path ( 0.4, -1.6) coordinate (N5);
      %
      % triangles
      %
      \draw [] (N0) -- (N1) -- (N2) -- cycle;
      \draw [] (N0) -- (N2) -- (N3) -- cycle;
      \draw [fill = gray!12] (N0) -- (N3) -- (N4) -- cycle;
      \draw [fill = gray!05] (N0) -- (N1) -- (N2) -- cycle;
      \draw [] (N0) -- (N4) -- (N5) -- cycle;
      \draw [] (N0) -- (N5) -- (N1) -- cycle;
      %
      % reference element
      \path ( -5.6, -1.2) coordinate (Q1);
      \path ( -5.6,  1.4) coordinate (Q2);
      \path ( -3.0, -1.2) coordinate (Q3);
      \draw [fill = gray!12] (Q1) -- (Q2) -- (Q3) -- cycle;
      \filldraw [black] (Q1)  circle (2.0pt);
      \filldraw [black] (Q2)  circle (2.0pt);
      \filldraw [black] (Q3)  circle (2.0pt);
      %
      % edges and labels
      %
      \path (-5.0, -0.1) coordinate (Khat);
      \path (-0.8,  0.4) coordinate (K);
      \node [below] at (Khat)  {$\hat{K}$};
      \node [below] at (K)  {$K$};
      \node [below] at (0.9, 0.8) {$K'$};
      \node [below right = 0pt and 0.3ex]  at (N0)  {$\bx_i$};
      %\node [below] at (0.2, -1.7)  {\small patch $\omega_i = \supp \phi_i$};
      \node [below] at (N5)  {\small patch $\omega_i = \supp \phi_i$};
      \filldraw [black] (N0)  circle (3.0pt);
      \filldraw [black] (N3)  circle (2.0pt);
      \filldraw [black] (N4)  circle (2.0pt);
      %
      % arrows
      %
      \path [->, line width = 0.75pt] ($(Khat) + (0.3,-0.2)$)
         edge [bend left] node [above] {$\bx = F_K(\bxi)$} ($(K) + (-0.3,-0.20)$);
   \end{tikzpicture}%
   \caption{Linear finite elements: reference element $\hat{K}$, mesh elements $K$ and $K'$, mapping $F_K$,
      node $\bx_i$, and the patch $\omega_i$ being the support of the
      corresponding basis function $\phi_i$.%
   }\label{fig:notation}
\end{figure}
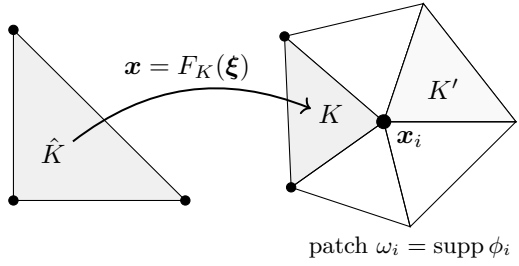
Later,  Kamenski, Huang, and Xu~\cite{KamHuaXu14} provided a bound for second-order elliptic PDEs for arbitrary conforming simplicial meshes without any conditions on the mesh regularity, which can be seen as a further generalization of the results above.
For $d \ge 1$ dimensions, their bound is
\begin{equation}
   \lambda_{\min}(A) \gtrsim
   \begin{cases}
   N^{-1},
   & d = 1,
   \\
   N^{-1} {\left(1 + \abs{\log \left( N \abs{K_{\min}} \right)} \right)}^{-1},
   & d = 2,
   \\
   {\left(  \sum\limits_{K\in\Th}  \Abs{K}^{1-\frac{d}{2}} \right)}^{-\frac{2}{d}},
   & d \ge 3.
   \end{cases}
   %\tag{KHX}
   \label{eq:kam}
\end{equation}
being some average of the mesh element volumes $\Abs{K}$~\cite[Lemma~5.1]{KamHuaXu14};
a comparable result was proven by Kamenski and Huang~\cite{KamHua14a} using a different approach.
For meshes fulfilling \cref{ass:grmcl}, \cref{eq:gramcl} and \cref{eq:kam} are comparable because the volume of a mesh element is comparable to the corresponding patch volume for such meshes.
Bound \cref{eq:kam} does not require any mesh regularity condition and, in contrast to \cref{eq:gramcl}, remains valid for arbitrary nonuniform meshes but it's elementwise averaging is not optimal.
For example, in the case of an internal layer with a thickness of only one element (\cref{fig:layer}), it holds that $\Abs{K_{\min}} \to 0$ but $\Abs{\omega_{\min}} > 1/2$ as the layer width goes to zero.

\begin{figure}[t]%
   \centering{}%
   \begin{tikzpicture}[scale = 1.37]%
      % define vertices 1, 2, 3, 4
      \tkzDefPoint(0, 0){p00}
      \tkzDefPoint(0, 1){p01}
      \tkzDefPoint(0, 2){p02}
      \tkzDefPoint(1, 0){p10}
      \tkzDefPoint(1, 1){p11}
      \tkzDefPoint(1, 2){p12}
      \tkzDefPoint(2, 0){p20}
      \tkzDefPoint(2, 1){p21}
      \tkzDefPoint(2, 2){p22}
      \tkzDefPoint(2.4, 0){p30}
      \tkzDefPoint(2.4, 1){p31}
      \tkzDefPoint(2.4, 2){p32}
      \tkzDefPoint(3.4, 0){p40}
      \tkzDefPoint(3.4, 1){p41}
      \tkzDefPoint(3.4, 2){p42}
      \tkzDefPoint(4.4, 0){p50}
      \tkzDefPoint(4.4, 1){p51}
      \tkzDefPoint(4.4, 2){p52}
      % description
      \coordinate (A) at (2.01, 2.2);
      \coordinate (B) at (2.41, 2.2);
      \coordinate (C) at (1.01, 2.2);
      \coordinate (D) at (4.6, 1.01);
      \draw [<->, line width = 0.75pt] (A) -- node [above] {$\varepsilon$} ++ (0.38,0);
      \draw [<->, line width = 0.75pt] (B) -- node [above] {$1$} ++ (0.98,0);
      \draw [<->, line width = 0.75pt] (C) -- node [above] {$1$} ++ (0.98,0);
      \draw [<->, line width = 0.75pt] (D) -- node [right] {$1$} ++ (0.0,0.98);
       % areas
      \tkzFillPolygon[gray!15](p20,p30,p32,p22)
      %% box + diagonals
      \tkzDrawSegments[thin](p00,p50 p50,p52 p52,p02 p02,p00) % box
      \tkzDrawSegments[thin](p01,p51 p02,p52 p10,p12 p20,p22 p30,p32 p30,p32 p40,p42) % grid
      \tkzDrawSegments[thin](p00,p11 p10,p21 p20,p31 p30,p41 p40,p51) % diagonals
      \tkzDrawSegments[thin](p01,p12 p11,p22 p21,p32 p31,p42 p41,p52) % diagonals
   \end{tikzpicture}%
   \caption{As $\varepsilon \to 0$, $\Abs{K_{\min}} = \varepsilon/2 \to 0$ 
      whereas $\Abs{\omega_{\min}} > 1/2$.%
   }\label{fig:layer}%
\end{figure}

Thus, both~\cite[Lemma~5.2]{GraMcL06} and~\cite[Lemma~5.1]{KamHuaXu14} are suboptimal.
The former has at least one additional factor depending on the mesh regularity whereas the latter has the suboptimal elementwise averaging.

The aim of this work is to show that \cref{ass:grmcl} is not necessary for obtaining \cref{eq:gramcl} (at least for the Lagrangian $P_m$ finite elements with simplicial meshes considered in this paper) and the factors depending on $F$, $G$, $H$, and $M$ can be completely removed from the bound \emph{without imposing any mesh regularity assumptions, neither global nor local}.

%%% Notation %%%%%%%%%%%%%%%%%%%%%%%%%%%%%%%%%%%%%%%%%%%%%%%%%%%%%%%%%%
\section{Considered problem, assumptions, and notation}\label{sec:notation}

Consider the boundary value problem (BVP) of a general diffusion differential equation
\begin{equation}
   \begin{cases}
      - \nabla \cdot \left( \D \nabla u \right) = f & \text{in $\Omega$}, \\
      u = 0                            & \text{on $\partial\Omega$},
   \end{cases}
   \label{eq:bvp1}
\end{equation}
where $\Omega \subset \mathbb{R}^d$ ($d \ge 2$) is a bounded Lipschitz domain.
It is assumed that the diffusion matrix $\D = \D (\bx)$ is symmetric and positive definite and there exist two constants $d_{\min}, d_{\max} > 0$ such that
\begin{equation}
   d_{\min} I \leq \D(\bx) \leq d_{\max} I \quad \forall \bx \in \Omega,
   \label{eq:D:1}
\end{equation}
where the $\ge$ sign means that the difference between the right-hand and left-hand side terms is positive semidefinite.

Let $\Th$ be a conforming, nondegenerate simplicial mesh (discretization) of $\Omega$.
For a mesh element $K \in \Th$, let $F_K \colon \hat{K} \to K$ be the invertible affine mapping from the reference mesh element $\hat{K}$ to $K$ (\cref{fig:notation}).
For simplicity, the reference element $\hat{K}$ is assumed to be unitary, i.e.,\ $\abs{\hat{K}} = 1$, so that the (constant) Jacobian matrix $F'_K$ satisfies $\Abs{\det(F_K')} = \Abs{K}$.

$V_h \subset H^1_0(\Omega)$ denotes the corresponding Lagrangian $\mathbb{P}_m$ ($m \ge 1$) finite element space with a given set of nodes
\[
   \set{\bx_i : i \in \cN},
   \quad \cN = \set{1, \dotsc, N},
\]
and the basis functions
\[
   \set{\phi_i : i \in \cN}
   ~~ \text{with} ~~
   \phi_i(\bx_j) = \delta_{ij}
\]
and
\[
   \supp \phi_i \subseteq \omega_i := \bigcup \set{ \bar{K}: \bx_i \in \bar{K} }
   ~~
   \forall i,j \in \cN
   .
\]
The finite element solution of the BVP~\cref{eq:bvp1} is defined by the variational formulation
\begin{equation}
   \int_\Omega \nabla v_h \cdot \D \nabla u_h \dx 
   = \int_\Omega f v_h \dx \quad \forall v_h \in V_h
   ,
   \label{eq:bvp:weak}
\end{equation}
which is equivalent to solving the system of linear equations
\[
   A \bu = \bm{f},
\]
where $A$ and $\bm{f}$ are defined via
\[
   A_{ij} := \int_\Omega \nabla\phi_i \cdot \D \nabla\phi_j \dx
   ~~\text{and}~~
   \bm{f}_i := \int_\Omega f \phi_i \dx \quad \forall i,j \in \cN
   .
\]
The solution of \cref{eq:bvp:weak} is then given by 
\[
   u_h =  \sum_{i \in \cN} \bu_i \phi_i
   .
\]
In the following, we are interested in the lower bound on the smallest eigenvalue $\lambda_{\min}(A)$ of $A$.

%%% d = 3 %%%%%%%%%%%%%%%%%%%%%%%%%%%%%%%%%%%%%%%%%%%%%%%%%%%%%%%%%%
\section{The general case \texorpdfstring{$d \ge 3$}{d greater than or equal to three}}
Hereafter, $a \gtrsim b$ means $a \ge C \cdot b$, where $C$ is a generic constant which can have different values at different appearances but is independent of the mesh and the solution of the BVP.\@

\begin{theorem}[$d \ge 3$]\label{thm:lmin}
Under the assumptions of \cref{sec:notation}, the smallest eigenvalue of the stiffness matrix $A$ for the finite element approximation of BVP~\cref{eq:bvp1} is bounded by
\begin{equation}
   \lambda_{\min}(A) \gtrsim
      {\left(  \sum\limits_{i\in\cN}  \Abs{\omega_i}^{1-\frac{d}{2}} \right)}^{-\frac{2}{d}}
      .
   \label{eq:lambdaMin}
\end{equation}
\end{theorem}

\begin{proof}
   The central idea for the proof is similar to~\cite[Lemma~5.1]{KamHuaXu14}: Poincaré and Sobolev inequalities, norm equivalence for the finite element functions over the reference element, and Hölder's inequality.
However, in order to arrive at the bound \cref{eq:lambdaMin}, some parts of the proof have to be done differently, namely, using a different norm on the reference element and doing this before applying Hölder's inequality.

First, assumption~\cref{eq:D:1} yields
\[
   \bu^T A \bu
      = \int_\Omega \nabla u_h \cdot \D \nabla u_h \dx
      \ge d_{\min} \Abs{u_h}^2_{H^1(\Omega)}
   .
\]
The next steps go back to Bank and Scott~\cite{BanSco89}: going from the $H^1$-seminorm through the $H^1$-norm (Poincaré inequality) to the $L^{2d/(d-2)}$-norm (Sobolev inequality),
\begin{align}
   \Abs{u_h}^2_{H^1(\Omega)}
   &\ge \frac{d_{\min} C_P}{1+C_P} \Norm{u_h}^2_{H^1(\Omega)}
   \notag{}
   \\
   &\ge \frac{d_{\min} C_P C_S}{1+C_P} \Norm{u_h}^2_{L^{\frac{2d}{d-2}}(\Omega)}
   .
   \label{eq:proofSobolev}
\end{align}
Then, represent the norm as an elementwise sum and use the transformation to the reference element $\hat{K}$,
\begin{align*}
   \Norm{u_h}^2_{L^{\frac{2d}{d-2}}(\Omega)}
   &= {\left( \int_\Omega \Abs{u_h}^{\frac{2d}{d-2}} \dx \right)}^{\frac{d-2}{d}}
   \\
   &= {\left( \sum_{K \in \Th} \int_K \Abs{u_h}^{\frac{2d}{d-2}} \dx \right)}^{\frac{d-2}{d}}
   \\
   &= {\left( \sum_{K \in \Th} \int_{\hat{K}} 
      \Abs{u_h \circ F_K}^{\frac{2d}{d-2}} \Abs{\det F'_K} \dxi 
      \right)}^{\frac{d-2}{d}}
   \\
   &= {\left( \sum_{K \in \Th} \Abs{K} \int_{\hat{K}}
      \Abs{u_h \circ F_K}^{\frac{2d}{d-2}}  \dxi \right)}^{\frac{d-2}{d}}
   \\
   &= {\left( \sum_{K \in \Th} 
      \Abs{K} \Norm{u_h \circ F_K}^{\frac{2d}{d-2}}_{L^{\frac{2d}{d-2}}(\hat{K})}
      \right)}^{\frac{d-2}{d}}
   .
\end{align*}

Norm equivalence for the finite-dimensional linear space of finite element functions over $\hat{K}$ yields
\[
   \Norm{u_h \circ F_K}^{\frac{2d}{d-2}}_{L^{\frac{2d}{d-2}}(\hat{K})}
   \gtrsim \Norm{\bm{u}_K}^{\frac{2d}{d-2}}_{l^{\frac{2d}{d-2}}}
   = \sum_{i \in \cN_K} \Abs{\bu_i}^{\frac{2d}{d-2}}
   ,
\]
where $\cN_K := \set{i : \bx_i \in \bar{K} }$ and $\bu_K := (\bu_i : i \in \cN_K)$ are the restrictions of $\cN$ and $\bu$ on $K$.
Hence,
\begin{align*}
   \bu^T A \bu
   &\gtrsim {\left( \sum_{K \in \Th}
      \Abs{K} \sum_{i \in \cN_K} \Abs{\bu_i}^{\frac{2d}{d-2}}
      \right)}^{\frac{d-2}{d}}
   \\
   &= {\left( \sum_{i \in \cN} \Abs{\bu_i}^{\frac{2d}{d-2}} \sum_{K \in \omega_i} \Abs{K}
      \right)}^{\frac{d-2}{d}}
   \\
   &= {\left( \sum_{i \in \cN} \Abs{\bu_i}^{\frac{2d}{d-2}} \Abs{\omega_i}
      \right)}^{\frac{d-2}{d}}
   \\
   &= {\left( \sum_{i \in \cN} {\left( \Abs{\bu_i}^2 \Abs{\omega_i}^\frac{d-2}{d} 
   \right)}^{\frac{d}{d-2}} \right)}^{\frac{d-2}{d}}
   .
\end{align*}
Applying Hölder's inequality for $1/p + 1/q = 1$ ($p, q > 1$) in the form
\[
   {\left(\sum_i \Abs{\alpha_i}^p\right)}^{\frac{1}{p}} \ge  \sum_i \Abs{\alpha_i \beta_i}
   {\left(\sum_i \Abs{\beta_i}^q\right)}^{-\frac{1}{q}}
   ,
   \quad \sum_i \Abs{\beta_i} \neq 0
\]
with
\[
   \alpha_i = \bu_i^2 \Abs{\omega_i}^\frac{d-2}{d},
   \quad \beta_i = \Abs{\omega_i}^{-\frac{d-2}{d}},
   \quad p = \frac{d}{d-2},
   \quad q = \frac{d}{2},
\]
finally gives
\begin{align*}
   \bu^T A \bu 
   &\gtrsim \sum_{i \in \cN} \Abs{\bu_i}^2
   {\left( \sum_{i \in \cN} {\left( \Abs{\omega_i}^{-\frac{d-2}{d}} \right)}^{\frac{d}{2}} \right)}^{-\frac{2}{d}}
   \\
   &= \sum_{i \in \cN} \Abs{\bu_i}^2 {\left( \sum_{i \in \cN} \Abs{\omega_i}^{1 -\frac{d}{2}} \right)}^{-\frac{2}{d}}
   .
   \qedhere
\end{align*}
\end{proof}

\begin{remark}[Geometric interpretation]
Bound \cref{eq:lambdaMin} can be rewritten as
\begin{align}
   \lambda_{\min}(A) \gtrsim N^{-1}
   \times
   {\left(
      \frac{1}{N} \sum_{i \in \cN}
      {\left( \frac{\Abs{\tilde\omega}}{\Abs{\omega_i}}\right)}^{\frac{d}{2} - 1}
   \right)}^{-\frac{2}{d}}
   ,
   \label{eq:geo:interpretation}
\end{align}
where $\Abs{\tilde\omega} := d\Abs{\Omega}/N$ can be seen as an equivalent to the average mesh patch size.
Alternatively, bound \cref{eq:geo:interpretation} also reads as

\[
   \lambda_{\min}(A)
   \gtrsim N^{-1}
   \times M_{1 - d/2}^{1 - 2/d}
      \left(\frac{\Abs{\tilde\omega}}{\Abs{\omega_i}}\right)
      ,
\]
where $M$ denotes the \emph{Hölder mean}
\[
   M_p (x_1,\dotsc, x_n) = {\left(\frac{1}{n} \sum\limits_{i = 1}^{n} x_i^p \right)}^{1/p}
   .
\]

Thus, the lower bound on the smallest eigenvalue is proportional to the inverse of the number of mesh vertices, $N^{-1}$, and to the power $1-2/d$ of the $d/2 - 1$ generalized mean patch volume irregularity.
\end{remark}

\begin{remark}%[Essential assumptions]
The essential assumptions for the proof are the following:
\begin{enumerate}[label=\alph*)]
\item symmetric, positive definite coefficient matrix $\D$ with $\lambda_{\min}(\D) \ge d_{\min} > 0$;
   \item validity of the Poincaré and Sobolev inequalities;
   \item conforming, nondegenerate simplicial meshes;
   no further mesh regularity conditions are required, neither global nor local.
\end{enumerate}
\end{remark}

%%% d = 2 %%%%%%%%%%%%%%%%%%%%%%%%%%%%%%%%%%%%%%%%%%%%%%%%%%%%%%%%%%
\section{The special case \texorpdfstring{$d=2$}{d equal to two}}

In two dimensions, the proof is mostly the same except that the Sobolev embedding $H^1(\Omega) \subset L^p(\Omega)$ holds for all finite $p$ and
\begin{equation}
   \Norm{u_h}^2_{H^1(\Omega)} \gtrsim \frac{1}{p} \Norm{u_h}^2_{L^p(\Omega)}
   \quad \forall u_h \in H^1(\Omega), ~ p < \infty.
   \label{eq:sobolev:2d}
\end{equation}

\begin{theorem}[$d = 2$]\label{thm:lmin:2d}
Under the assumptions of \cref{sec:notation}, the smallest eigenvalue of the stiffness matrix $A$ for the finite element approximation of BVP~\cref{eq:bvp1} is bounded by
\begin{equation}
   \lambda_{\min}(A) 
      \gtrsim N^{-1} {\left( 1 + \bigAbs{\ln (N \Abs{\omega_{\min}})} \right)}^{-1}.
    \label{eq:lambdaMin:2d}
\end{equation}
\end{theorem}

\begin{proof}
Repeating the steps of the proof of \cref{thm:lmin} and using \cref{eq:sobolev:2d} instead of \cref{eq:proofSobolev} yield for $2 < p < \infty$
\allowdisplaybreaks{}%
\begin{align}
   \bu^T A \bu
   &\gtrsim \frac{1}{p} \Norm{u_h}^2_{L^p(\Omega)}
   \tag{Poincaré, Sobolev}
   \\
   &= \frac{1}{p} {\left( \sum_{K \in \Th} \Norm{u_h}^p_{L^p(K)}\right)}^{\frac{2}{p}}
  \notag
   \\
   &= \frac{1}{p} {\left( \sum_{K \in \Th} \Abs{K} \Norm{u_h \circ F_K}^p_{L^p(\hat{K})}\right)}^{\frac{2}{p}}
   \notag
   \\
   &\gtrsim \frac{1}{p} {\left( \sum_{K \in \Th} \Abs{K} \Norm{\bu_K}^p_{l^p}\right)}^{\frac{2}{p}}
   \tag{norm equivalence}
   \\
   &=  \frac{1}{p} {\left( \sum_{i \in \cN} \bu_i^p \Abs{\omega_i} \right)}^{\frac{2}{p}}
   \notag
   \\
   &\ge \frac{1}{p} \sum_{i \in \cN} \Abs{\bu_i}^2
      {\left( \sum_{i \in \cN} \Abs{\omega_i}^{-\frac{2}{p-2}} \right)}^{-\frac{p-2}{p}}
   \text{(Hölder's)}
   \label{eq:before:rough:estimate}
   \\
   &\ge \sum_{i \in \cN} \Abs{\bu_i}^2 \frac{1}{p} {\left( N \Abs{\omega_{\min}}^{-\frac{2}{p-2}} \right)}^{-\frac{p-2}{p}}
   \label{eq:rough:estimate}
   \\
   &= \sum_{i \in \cN} \Abs{\bu_i}^2 N^{-1} \frac{1}{p} {\left( N \Abs{\omega_{\min}} \right)}^{\frac{2}{p}}
   \notag
   .
\end{align}
Choosing $p = \max\set{2, \Abs{\log (N \Abs{\omega_{\min}})}}$ provides the largest lower bound
\[
   \bu^T A \bu
   \gtrsim 
      \sum_{i \in \cN} \Abs{\bu_i}^2 N^{-1} 
      {\left( 1 + \bigAbs{\ln (N \Abs{\omega_{\min}})} \right)}^{-1}
   .
\]
For the choice $p=2$, the bound is defined as the limiting case for $p \to 2^+$.
\end{proof}

\begin{remark}[a practical bound in two dimensions]\label{rem:2d}
In addition to the factor $N^{-1}$, bound \cref{eq:lambdaMin:2d} also involves the logarithmic factor $\log (N \Abs{\omega_{\min}})$.
In many practical applications, neither $N$ nor $\Abs{\omega_{\min}}$ can become arbitrarily large or small since it is not reasonable to construct meshes with more than \num{e10} elements or to refine below \num{e-10} (the size of a single atom).
Thus, for the most two-dimensional application scenarios, we can assume $\log (N \Abs{\omega_{\min}}) = \mathcal{O}(1)$ and consider the bound as $\lambda_{\min}(A) \gtrsim N^{-1}$.%
\footnote{Except for the unfortunate cases when meshing software produce degenerate elements of vanishing size.}
This is also confirmed by numerical experiments in \cref{sec:examples}: for all considered two-dimensional examples the behaviour of the exactly computed $\lambda_{\min}(A)$ is almost indistinguishable from $N^{-1}$ (\cref{fig:uniform:2d,fig:n:shishkin,fig:n:shishkin2,fig:n:bakhvalov,fig:n:bakhvalov2,fig:n:power,fig:n:single}).
\end{remark}
%%%% Numerical examples %%%%%%%%%%%%%%%%%%%%%%%%%%%%%%%%%%%%%%%%%%%%%%%
\section{Numerical experiments}\label{sec:examples}

%---Calibration -----------------------------------------------
\begin{figure*}[t]
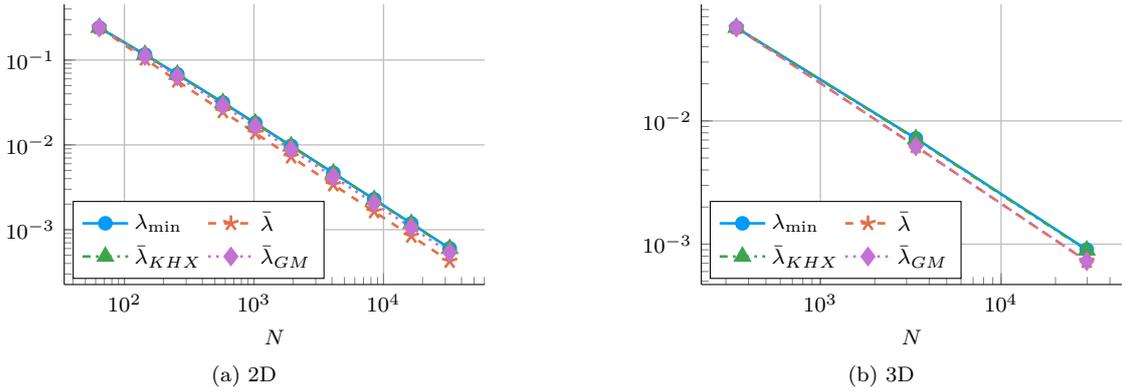
%
  \centering{}%
  \subfloat[2D]{%
    \PlotResults{n_uniform.dat}{$N$}{0.33}%
    \label{fig:uniform:2d}%
  }%
  \hspace{0.1\linewidth}%
  \subfloat[3D]{%
    \PlotResults{n_uniform_3d.dat}{$N$}{0.33}%
  }%
  \caption{Constants are chosen such that the bounds coincide with the exact value on uniform meshes.}%
  \label{fig:n:uniform}%
\end{figure*}

In the following, we compare the estimates of~\cref{thm:lmin,thm:lmin:2d} with \cref{eq:gramcl2,eq:kam} by computing the estimates $\bar{\lambda}$, $\bar{\lambda}_{GM}$, and $\bar{\lambda}_{KHX}$ of \cref{thm:lmin:2d}, \cite[Lemma~5.2]{GraMcL06}, and~\cite[Lemma~5.1]{KamHuaXu14}, respectively,
\begin{subequations}
\begin{align}
  \bar{\lambda} &\approx
  N^{-1} {\left( 1 + \bigAbs{\ln (N \Abs{\omega_{\min}})} \right)}^{-1}
  ,
  \label{eq:k2}
  \\
  \bar{\lambda}_{GM} &\approx
  N^{-1} {\left( 1 + \bigAbs{\ln \big({(MH)}^{-1} N \Abs{\omega_{\min}} \big)} \right)}^{-1}
  ,
  \label{eq:gm2}
  \\
  \bar{\lambda}_{KHX} &\approx
  N^{-1}_{ele} {\left(1 + \abs{\log \left( N_{ele} \abs{K_{\min}} \right)} \right)}^{-1}
  \label{eq:khx2}
\end{align}
\end{subequations}
in two dimensions ($d = 2$) and
\begin{subequations}
\begin{align}
  \bar{\lambda} &\approx
  {\left(  \sum\limits_{i\in\cN}  \Abs{\omega_i}^{-\frac{1}{2}} \right)}^{-\frac{2}{3}}
  ,
  \label{eq:k3}
  \\
  \bar{\lambda}_{GM} &\approx
  \frac{1}{\sqrt[3]{M H}}
  {\left(  \sum\limits_{i\in\cN}  \Abs{\omega_i}^{-\frac{1}{2}} \right)}^{-\frac{2}{3}}
  ,
  \label{eq:gm3}
  \\
  \bar{\lambda}_{KHX} &\approx
  {\left(  \sum\limits_{K\in\Th}  \Abs{K}^{-\frac{1}{2}} \right)}^{-\frac{2}{3}}
  \label{eq:khx3}
\end{align}
\end{subequations}
in three dimensions ($d = 3$).
The (unknown) constants have been chosen empirically such that the corresponding estimates coincide with the exact $\lambda_{\min}$ in the case of uniform meshes (\cref{fig:n:uniform}).
Numerical computations were carried out with \emph{Julia}~v1.5.2~\cite{julia}.

Since the formulas for $\bar{\lambda}$ and $\bar{\lambda}_{GN}$ differ only by the factor involving $M$ and $H$, the difference in the numerical results for these two estimates also indicates the value of this factor.
In this sense, the essential difference between~\cite[Lemma~5.2]{GraMcL06} and \cref{thm:lmin} (and~\ref{thm:lmin:2d}) is that the former provides a bound on $\lambda_{\min}(A)$ which depends on $H$ whereas the latter provides a sharp bound which does not depend on $H$ or any other mesh regularity constants.

\begin{figure*}[p]
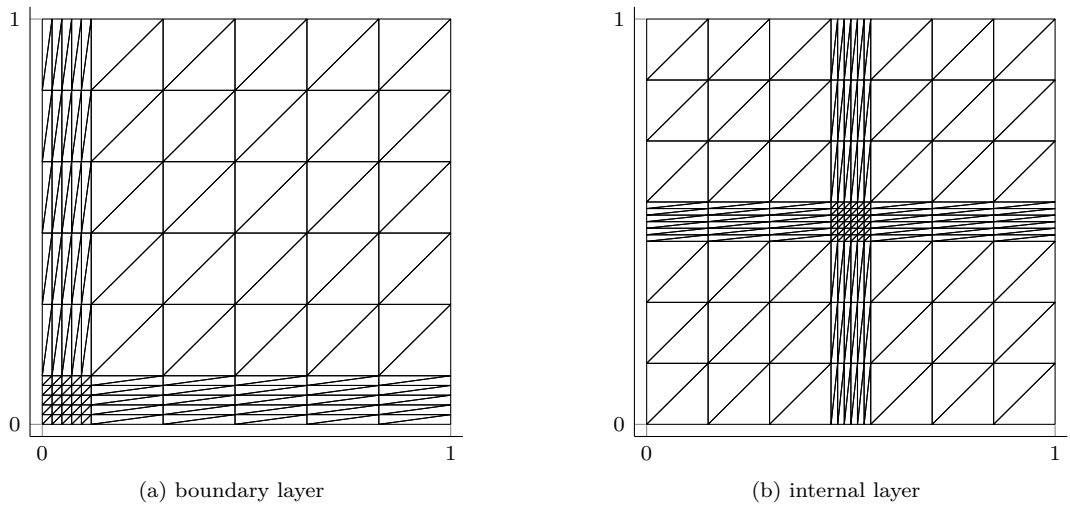

  \centering{}%
  \subfloat[boundary layer]{%
    \PlotMesh{shishkin}{0.36}%
    \label{fig:mesh:shishkin}%
  }%
  \hspace{0.1\linewidth}%
  \subfloat[internal layer]{%
    \PlotMesh{shishkin_2}{0.36}%
    \label{fig:mesh:shishkin2}%
  }%
  \caption{Shishkin meshes, $\varepsilon = 0.05$ (2D).}%
  \label{fig:mesh:s}
\end{figure*}

\begin{figure*}[p]
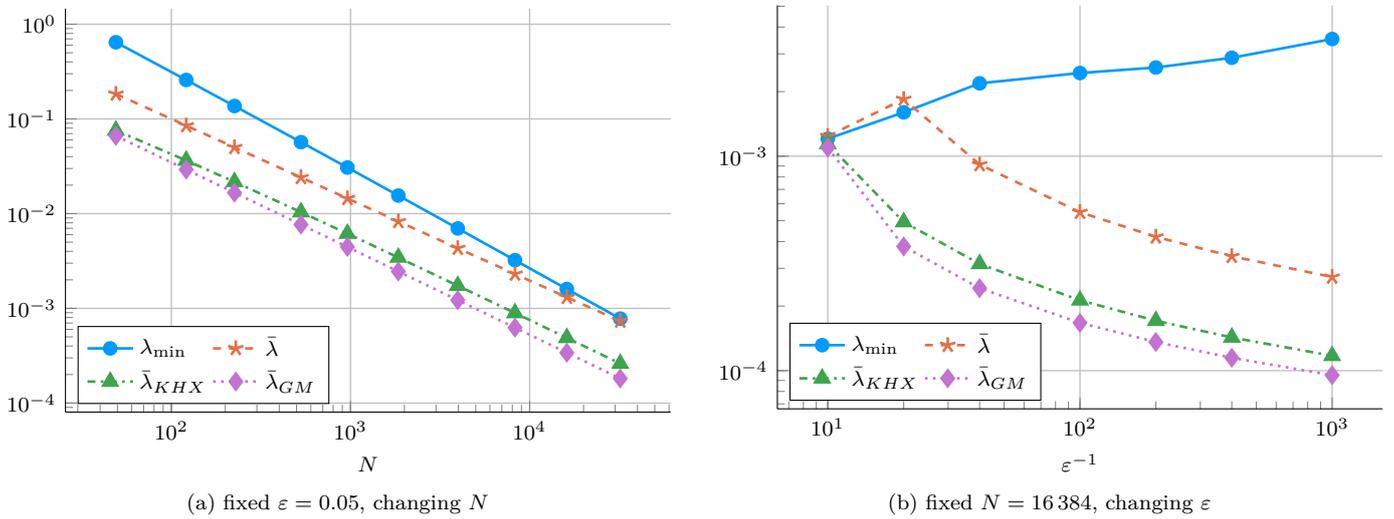

  \centering{}%
  \subfloat[fixed $\varepsilon = 0.05$, changing $N$]{%
    \PlotResults{n_shishkin.dat}{$N$}{0.475}%
    \label{fig:n:shishkin}%
  }%
  \hfill{}%
  \subfloat[fixed $N = \num{16384}$, changing $\varepsilon$]{%
    \PlotResults{e_shishkin.dat}{$\varepsilon^{-1}$}{0.475}%
    \label{fig:e:shishkin}%
  }%
  \caption{Shishkin meshes (2D boundary layer).}\label{fig:ne:shishkin}
\end{figure*}

\begin{figure*}[p]
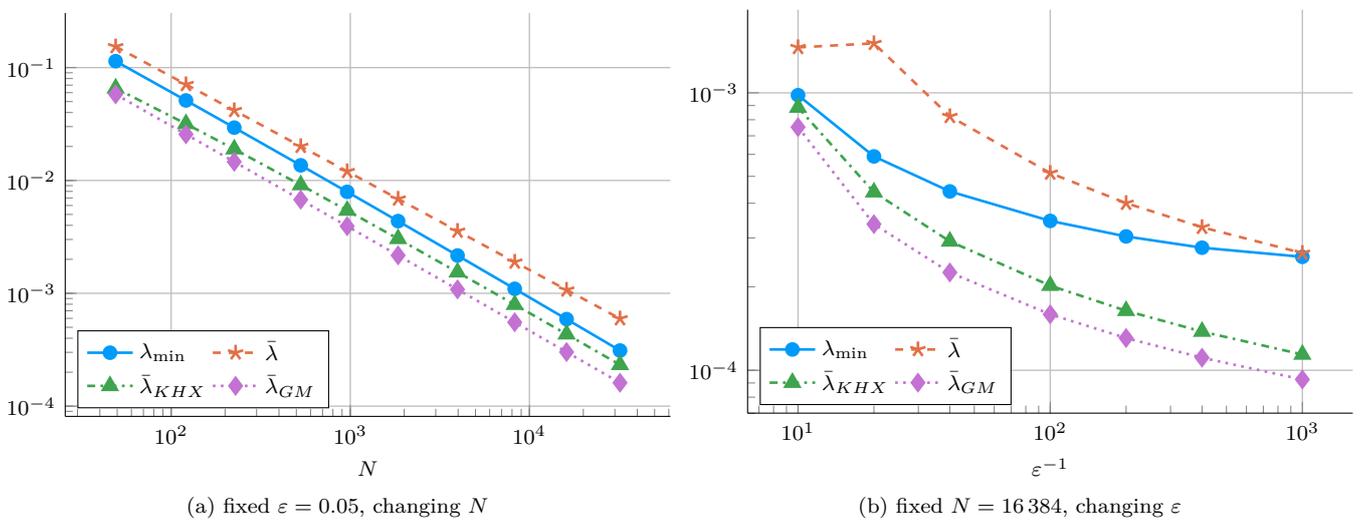

  \centering{}%
  \subfloat[fixed $\varepsilon = 0.05$, changing $N$]{%
    \PlotResults{n_shishkin_2.dat}{$N$}{0.475}%
    \label{fig:n:shishkin2}%
  }%
  \subfloat[fixed $N = \num{16384}$, changing $\varepsilon$]{%
    \PlotResults{e_shishkin_2.dat}{$\varepsilon^{-1}$}{0.475}%
    \label{fig:e:shishkin2}%
  }%
  \caption{Shishkin meshes (2D internal layer).}\label{fig:ne:shishkin2}
\end{figure*}

%--- Shishkin -----------------------------------------------
\begin{meshexample}[Shishkin]
   A one-dimensional Shishkin mesh is a layer-adapted mesh consisting of two equidistant parts: a  fine part and a coarse part (\cref{fig:mesh:s}).
   The nodes of the fine, boundary layer part are given by
\[
   x_i = \min\{1, 2 C_\sigma \varepsilon \ln N \} \times \frac{i}{N},
   \quad 0 \le i \le \frac{N}{2}
   .
\]
In the coarse part, the remaining nodes are spaced equally.
This mesh is not locally uniform and has an abrupt change in the mesh step size at the transition point.
A two dimensional Shishkin mesh (\cref{fig:mesh:shishkin}) is a product mesh
\[
   \boldsymbol{x}_{i,j} := (x_i, x_j), \qquad 0 \le i,j \le N
   .
\]

The step size change from the coarse part to the fine part of the one-dimensional Shishkin mesh is proportional to $\varepsilon \ln N$, which yields $H \sim \varepsilon^{-2} / \ln^2 N$ in two and $H \sim \varepsilon^{-3} / \ln^3 N$ in three dimensions.

\Cref{fig:ne:shishkin,fig:ne:shishkin2} show the results of the numerical experiment.
First, the mesh parameter $\varepsilon$ is fixed and the number of nodes $N$ is changing to verify the dependence of the estimate on $N$ (\cref{fig:n:shishkin,fig:n:shishkin2}).
Then, $N$ is fixed and the mesh grading parameter $\varepsilon$ is changed to observe the dependence of the estimate on the mesh shape and the change of the mesh step size (\cref{fig:e:shishkin,fig:e:shishkin2}).

All three estimates have similar behaviour but $\bar{\lambda}$ provides a noticeably larger bound (by a factor of approximately 3) than $\bar{\lambda}_{GM}$ and $\bar{\lambda}_{KHX}$.
Estimates $\bar{\lambda}_{GM}$ and $\bar{\lambda}_{KHX}$ are very similar for this case, $\bar{\lambda}_{KHX}$ being slightly better than $\bar{\lambda}_{GM}$.
\end{meshexample}

\begin{figure*}[p]
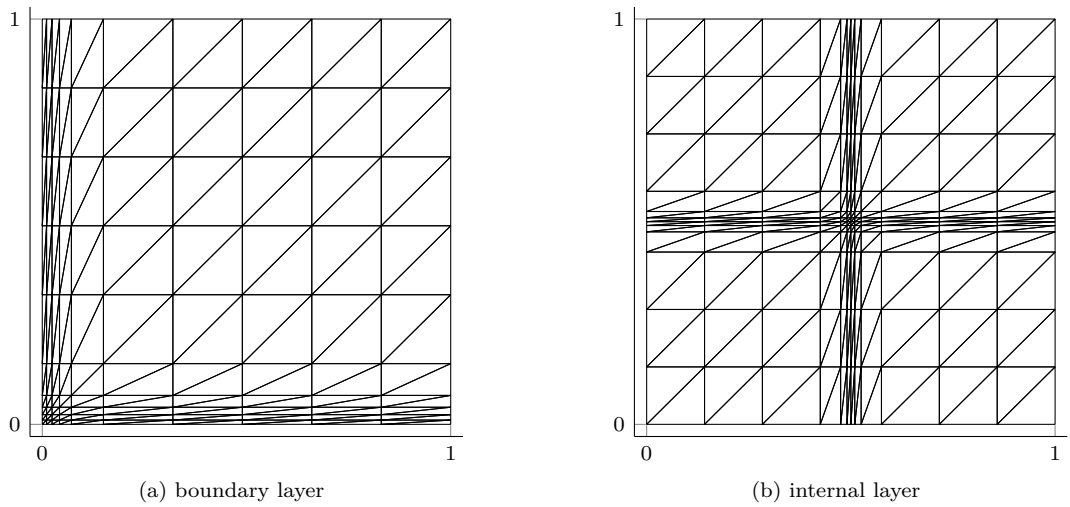

  \centering{}%
  \subfloat[boundary layer]{%
    \PlotMesh{bakhvalov_type}{0.36}%
    \label{fig:mesh:bakhvalov}%
  }%
  \hspace{0.1\linewidth}%
  \subfloat[internal layer]{%
    \PlotMesh{bakhvalov_type_2}{0.36}%
    \label{fig:mesh:bakhvalov2}%
  }%
  \caption{Bakhvalov-type meshes, $\varepsilon = 0.05$ (2D).}%
  \label{fig:mesh:b}%
\end{figure*}

\begin{figure*}[p]
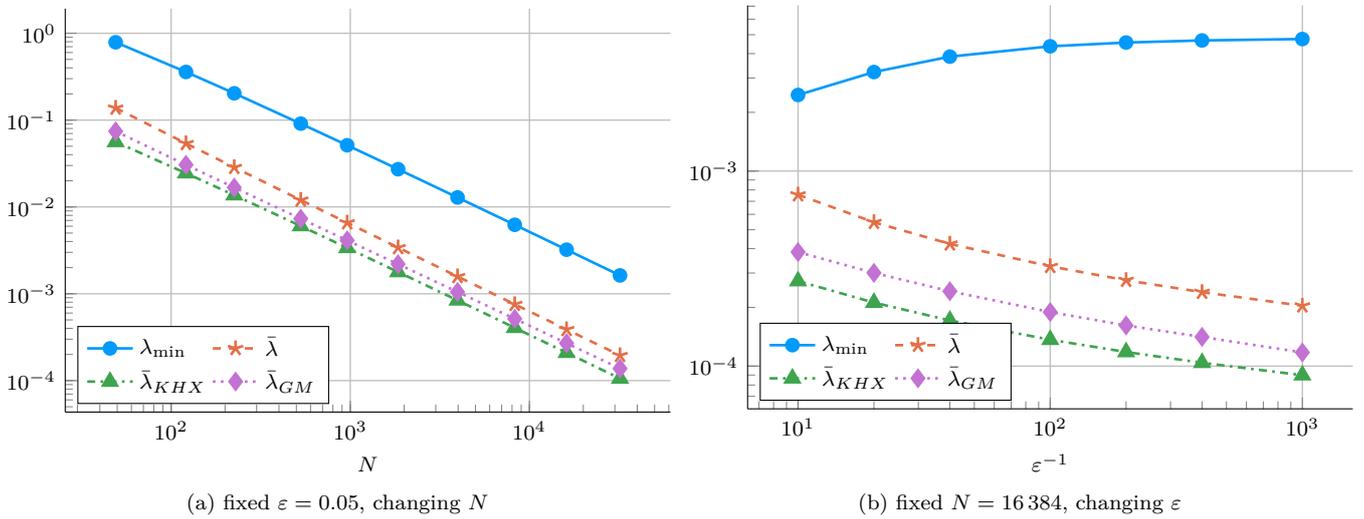

  \centering{}%
  \subfloat[fixed $\varepsilon = 0.05$, changing $N$]{%
    \PlotResults{n_bakhvalov_type.dat}{$N$}{0.475}%
    \label{fig:n:bakhvalov}%
  }%
  \subfloat[fixed $N = \num{16384}$, changing $\varepsilon$]{%
    \PlotResults{e_bakhvalov_type.dat}{$\varepsilon^{-1}$}{0.475}%
    \label{fig:e:bakhvalov}%
  }%
  \caption{Bakhvalov-type meshes (2D boundary layers).}\label{fig:ne:bakhvalov}
\end{figure*}

\begin{figure*}[p]
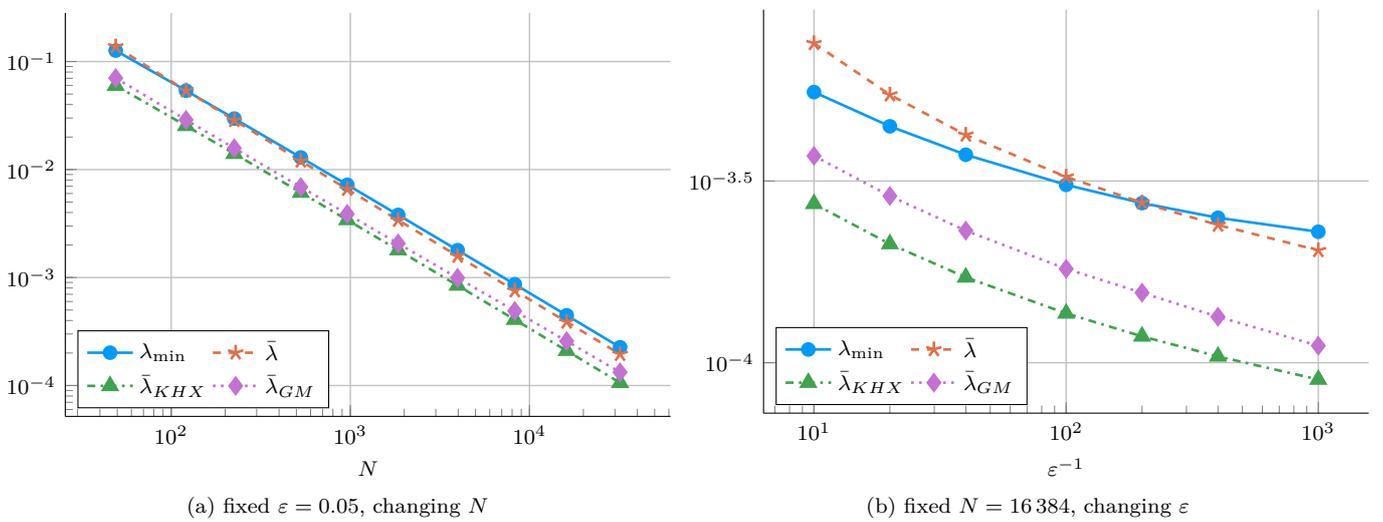

  \centering{}%
  \subfloat[fixed $\varepsilon = 0.05$, changing $N$]{%
    \PlotResults{n_bakhvalov_type_2.dat}{$N$}{0.475}%
    \label{fig:n:bakhvalov2}%
  }%
  \subfloat[fixed $N = \num{16384}$, changing $\varepsilon$]{%
    \PlotResults{e_bakhvalov_type_2.dat}{$\varepsilon^{-1}$}{0.475}%
    \label{fig:e:bakhvalov2}%
  }%
  \caption{Bakhvalov-type meshes (2D internal layers).}\label{fig:ne:bakhvalov2}
\end{figure*}

%--- Bakhvalov -----------------------------------------------
\begin{meshexample}[Bakhvalov type]\label{ex:bakhvalov}
A Bakhvalov type mesh (\cref{fig:mesh:b}) is defined by
\[
   x_i = -C_\sigma \varepsilon 
    \times \ln \left(1 - 2(1-\varepsilon)\frac{i}{N} \right),
   \qquad 0 \le i \le \frac{N}{2}
   ,
\]
in the fine, boundary layer part and is equidistant in the remaining coarse part.
The transition point from the fine to the coarse mesh is
$\sigma = \min\{0.5, - C_\sigma \varepsilon \ln \left(\varepsilon\right) \}$.
In the numerical tests, $C_\sigma = 1$ was used for both Shishkin and Bakhvalov type mesh examples.

Numerical comparison of the Bakhvalov type meshes (\cref{fig:ne:bakhvalov,fig:ne:bakhvalov2}) show similar results: $\bar{\lambda}$ provides a larger bound than $\bar{\lambda}_{GM}$ and $\bar{\lambda}_{KHX}$ by a factor of roughly 2, estimate $\bar{\lambda}_{GM}$ being better than $\bar{\lambda}_{KHX}$.
The less abrupt step size change of the Bakhvalov type mesh yields a smaller $H$ and, thus, a lesser underestimation for $\bar{\lambda}_{GM}$.
\end{meshexample}

\begin{figure*}[p]
  \centering{}%
  \subfloat[power-graded, $\beta = 3$]{%
     \includegraphics[width=0.35\linewidth]{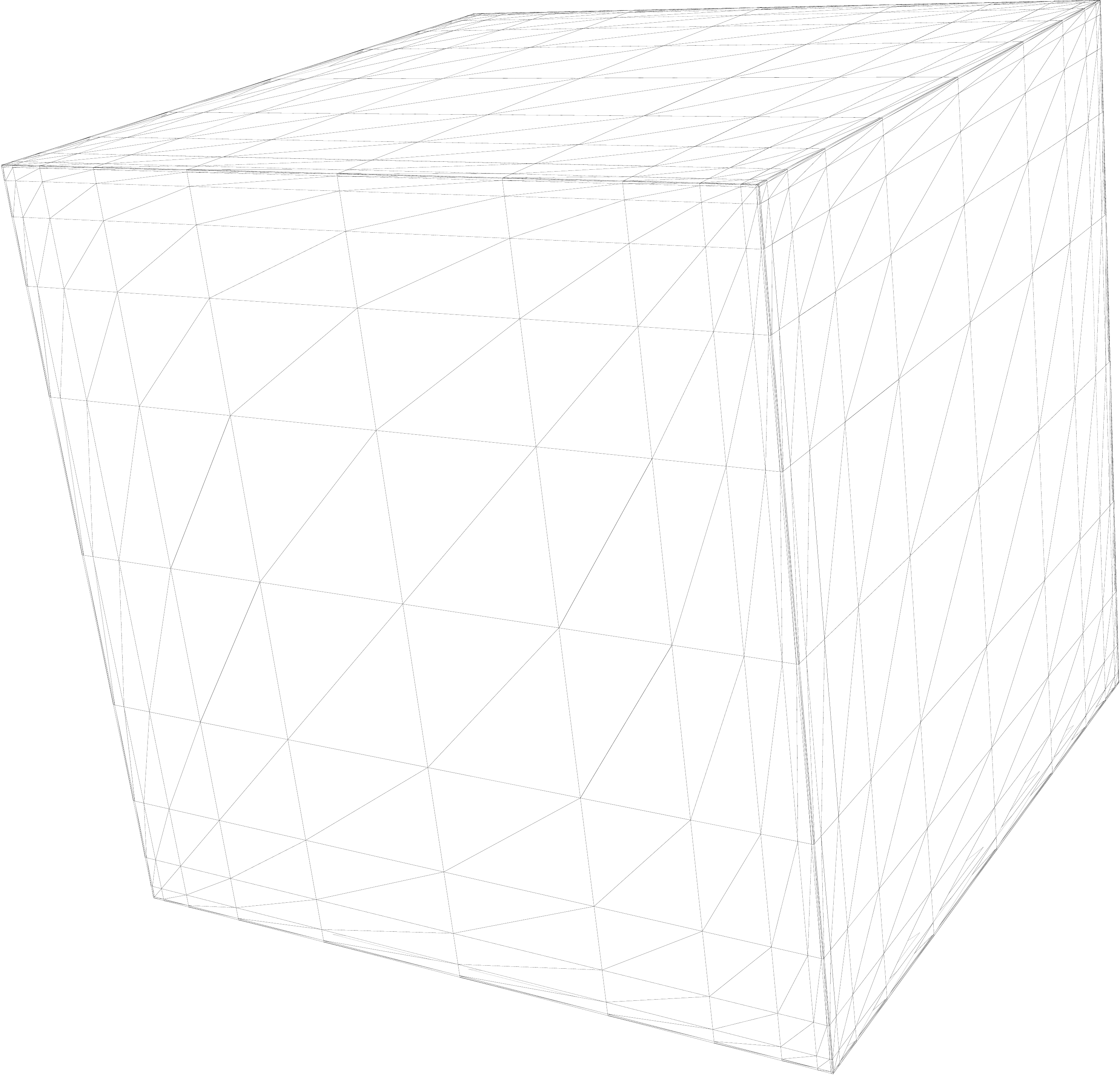}%
    \label{fig:mesh:power3}%
  }%
  \hspace{0.1\linewidth}%
  \subfloat[single element layer, $\varepsilon = 0.1$]{%
     \includegraphics[width=0.35\linewidth]{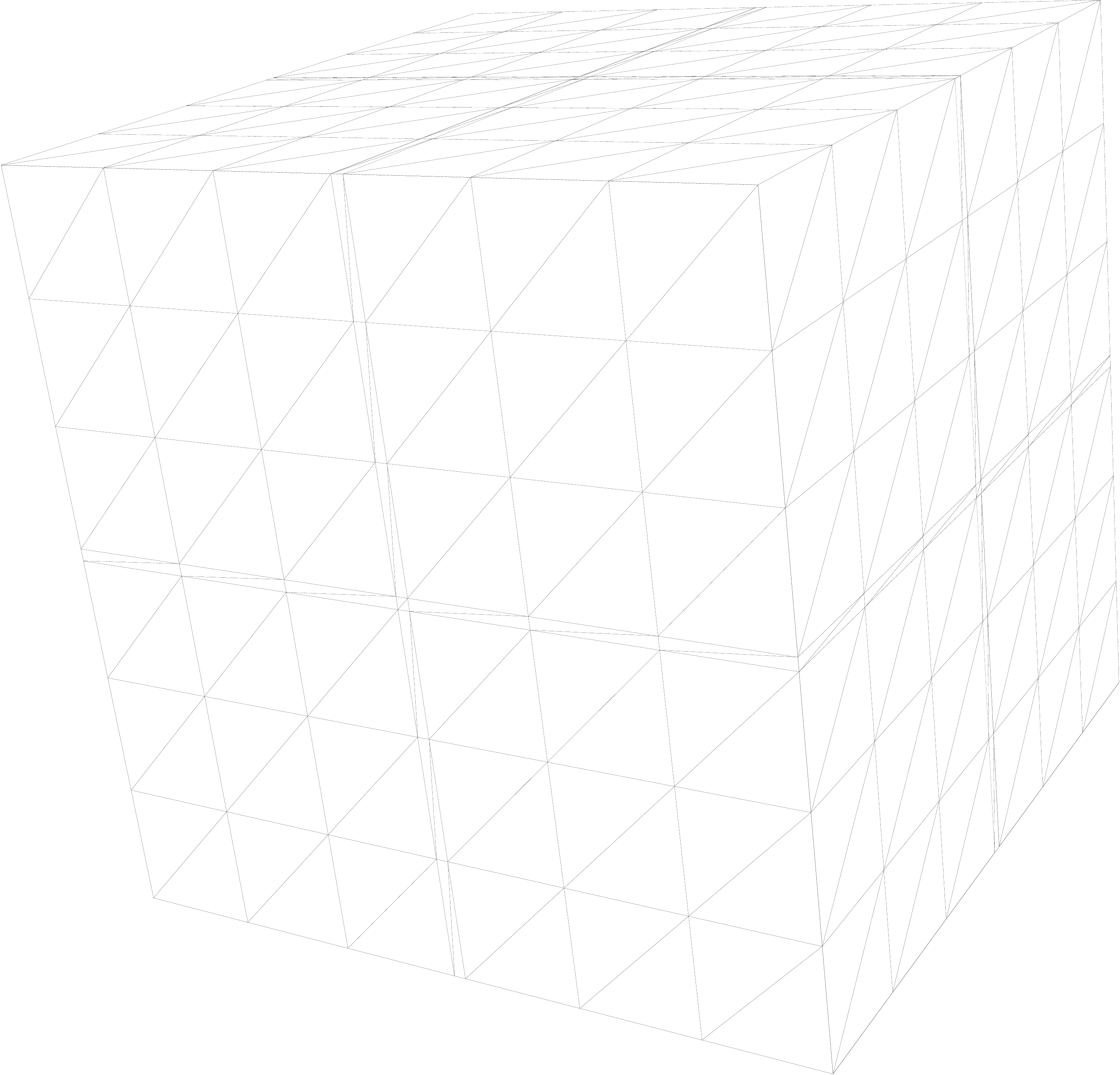}%
    \label{fig:mesh:single3}%
  }%
  \caption{Power-graded and single element internal layer meshes (2D).}%
  \label{fig:mesh:ps3}
\end{figure*}

%--- Power-graded -----------------------------------------------
\begin{meshexample}[power-graded]\label{ex:power}
   A power graded mesh (\cref{fig:mesh:power,fig:mesh:power3}) is defined by
\[
  x_i = \begin{cases}
    \frac{1}{2} \times {\left(2 \frac{i}{N} \right)}^\beta,
    & 0 \le i \le \frac{N}{2},
    \\
    1 - x_{N - i},
    & \frac{N}{2} < i \le N
    .
   \end{cases}
\]

Numerical results for the two-dimensional power-graded meshes (\cref{fig:ne:power}) again show that $\bar{\lambda}$ provides a larger bound than $\bar{\lambda}_{GM}$ and $\bar{\lambda}_{KHX}$ by a factor of roughly \numrange{1.5}{2.5}, $\bar{\lambda}_{GM}$ being better than $\bar{\lambda}_{KHX}$.

In three dimensions, the differences must be much more striking because the effect of the mesh nonuniformity is not eliminated by the logarithm in \cref{eq:k2,eq:gm2,eq:khx2} but appears in the estimates directly:  only elementwise averaging in the case of  $\bar{\lambda}_{KHX}$ in \cref{eq:k3} and a factor proportional to $1/\sqrt[3]{MH}$ for $\bar{\lambda}_{GM}$ in \cref{eq:gm3}.
\Cref{fig:ne:power3} confirms this consideration: $\bar{\lambda}$ provides a much better bound than both $\bar{\lambda}_{KHX}$ and $\bar{\lambda}_{GM}$ and, in the case of a fixed $N$ and changing $\beta$, it also reproduces the shape of the curve for the exact value more accurately.
In particular, it can be observed that the value $H$ can become very large for highly irregular anisotropic meshes but it does not affect the value of the smallest eigenvalue directly, as predicted by \cref{thm:lmin}.
\end{meshexample}

%---Single layer -----------------------------------------------
\begin{meshexample}[Single element layer]\label{ex:single}
This mesh simulates an internal layer with a thickness of only one element (\cref{fig:layer,fig:mesh:single,fig:mesh:single3}).
The remaining part is equidistant.
This example is considered for illustrative purposes to show the effect of a very abrupt change of the mesh step size.

Numerical results in \Cref{fig:ne:single,fig:ne:single3} show that $\bar{\lambda}$ is very close to the exact value, whereas $\bar{\lambda}_{GM}$ and $\bar{\lambda}_{KHX}$ underestimate the exact value by a large factor, which is proportional to $\varepsilon$ (due to the special construction of this example).
Moreover, it can be observed that a very abrupt step size change doesn't affect the smallest eigenvalue as long as only single elements are affected and not the whole patches (though it will affect the largest eigenvalue).
\end{meshexample}

\begin{remark}[Possible improvements]\label{rem:bl}
Results in \cref{fig:e:shishkin,fig:e:bakhvalov,fig:e:power} demonstrate an interesting phenomenon: $\lambda_{\min}$ increases slightly as $\varepsilon$ is decreasing and the mesh becomes more anisotropic.
The reason for this might be that, as $\varepsilon$ decreases, the mesh boundary layer becomes thinner and moves closer to the boundary.
It has been observed~\cite{KamHua14a} that mesh elements near the boundary have a lesser impact on the smallest eigenvalue than elements inside the domain, so that `pressing' the boundary layer containing a significant number of mesh elements towards the boundary might have a similar effect on $\lambda_{\min}$ as mesh coarsening.
Bounds of \cref{thm:lmin,thm:lmin:2d} do not take this effect into account, so that there still might be some room for further improvement for this particular case, although it is not clear if such bounds can have the same simple form.

Other than that, the new bounds are clearly sharper and seem to more or less follow the curve for the exact value (\cref{fig:ne:shishkin2,fig:ne:bakhvalov2,fig:ne:single}).

In three dimensions, it doesn't seem that the bound can be further improved in general, at least not in its current form, since its derivation is free from any mesh reguarity assumptions and rough estimates.

In two dimensions, the only rough estimate in the derivation happens bewtween \cref{eq:before:rough:estimate,eq:rough:estimate}.
Hence, a possible improvement might be a replacement of the minimal patch volume $\abs{\omega_{\min}}$ in the estimate with some average, following the structure of the bound for $d\ge3$.
However, it not clear to the author yet how to deal with the sum in \cref{eq:before:rough:estimate} directly or if it can be improved at all.
Considering that $\abs{\omega_{\min}}$ enters the estimate only inside the logarithm, this possible improvement would be only a small one.
In particular, as already predicted in \cref{rem:2d}, in all two-dimensional exmaples the curve for the exact smallest eigenvalue is comparable to $N^{-1}$ as $N$ increases.
\end{remark}

%--- Shishkin meshes -------------------------------------------------


%--- Bakhvalov-type meshes -------------------------------------------
%
%--- Power-graded and single layer -----------------------------------
\begin{figure*}[p]
  \centering{}%
  \subfloat[power-graded, $\beta = 3$]{%
    \PlotMesh{power_graded}{0.36}%
    \label{fig:mesh:power}%
  }%
  \hspace{0.1\linewidth}%
  \subfloat[single element layer, $\varepsilon = 0.2$]{%
    \PlotMesh{single_layer}{0.36}%
    \label{fig:mesh:single}%
  }%
  \caption{Power-graded and single element internal layer meshes (2D).}%
  \label{fig:mesh:ps}%
\end{figure*}
\begin{figure*}[p]
  \centering{}%
  \subfloat[fixed $\beta=3.0$, changing $N$]{%
    \PlotResults{n_power_graded.dat}{$N$}{0.475}%
    \label{fig:n:power}%
  }%
  \subfloat[fixed $N = \num{16384}$, changing $\beta$]{%
    \PlotResults{e_power_graded.dat}{$\beta$}{0.475}%
    \label{fig:e:power}%
  }%
  \caption{Power-graded meshes (2D).}\label{fig:ne:power}
\end{figure*}
\begin{figure*}[p]
  \centering{}%
  \subfloat[fixed $\varepsilon = 0.1$, changing $N$]{%
    \PlotResults{n_single_layer.dat}{$N$}{0.475}%
    \label{fig:n:single}%
  }%
  \subfloat[fixed $N = \num{16384}$, changing $\varepsilon$]{%
    \PlotResults{e_single_layer.dat}{$\varepsilon^{-1}$}{0.475}%
    \label{fig:e:single}%
  }%
  \caption{Single element internal layer (2D).}\label{fig:ne:single}%
\end{figure*}

%--- Convergence plots 3D -----------------------------------------------
\begin{figure*}[p]
  \centering{}%
  \subfloat[fixed $\beta=3.0$, changing $N$]{%
    \PlotResults{n_power_graded_3d.dat}{$N$}{0.475}%
    \label{fig:n:power3}%
  }%
  \subfloat[fixed $N = \num{1728}$, changing $\beta$]{%
    \PlotResults{e_power_graded_3d.dat}{$\beta$}{0.475}%
    \label{fig:e:power3}%
  }%
  \caption{Power-graded meshes (3D).}\label{fig:ne:power3}
\end{figure*}
\begin{figure*}[p]
  \centering{}%
  \subfloat[fixed $\varepsilon = 0.05$, changing $N$]{%
    \PlotResults{n_single_layer_3d.dat}{$N$}{0.475}%
    \label{fig:n:single3}%
  }%
  \subfloat[fixed $N = \num{1728}$, changing $\varepsilon$]{%
    \PlotResults{e_single_layer_3d.dat}{$\varepsilon^{-1}$}{0.475}%
    \label{fig:e:single3}%
  }%
  \caption{Single element internal layer (3D).}\label{fig:ne:single3}
\end{figure*}

\section*{Acknowledgments}
I would like to thank Mark Ainsworth for pointing out the work of Graham and McLean~\cite{GraMcL06} and sharing with me his point of view about their `unhappiness' with the local mesh regularity assumption in the proof.
This gave me an idea of how the bound should look and an additional motivation to find a way to remove this assumption.
Also, I would like to acknowledge the anonymous referee for the comments and very helpful suggestions which helped significantly to improve the quality of this paper.

%--- References -------------------------------------------------------
%\bibliographystyle{abbrv}
%\bibliography{kamenski-lmin}

\end{document}